\RequirePackage{fix-cm}
\documentclass[smallextended]{svjour3}       
\smartqed  
\usepackage{graphicx}
\usepackage{amsmath}
\usepackage{amssymb}
\usepackage{url}
\usepackage{color}
\usepackage[normalem]{ulem}

\begin{document}

\titlerunning{Curves that lie on a geodesic sphere or on a totally geodesic hypersurface}
\title{Characterization of curves that lie on a geodesic sphere or on a totally geodesic hypersurface in a hyperbolic space or in a sphere\thanks{This is a pre-print of an article published in Mediterranean Journal of Mathematics: \\
da Silva, L. C. B. $\&$ da Silva, J. D. Mediterr. J. Math. (2018) \textbf{15}: 70. https://doi.org/10.1007/s00009-018-1109-9}
}


\author{Luiz C. B. da Silva         \and
        Jos\'e Deibsom da Silva 
}


\institute{Da Silva, L. C. B.  \at
              Departamento de Matem\'atica, Universidade Federal de Pernambuco \\
             50670-901, Recife, Pernambuco, Brazil\\
              \email{luizsilva@dmat.ufpe.br}           
           \and
           Deibsom da Silva, J. \at
             Departamento de Matem\'atica, Universidade Federal Rural de Pernambuco, \\
         52171-900, Recife, Pernambuco, Brazil\\
         \email{jose.dsilva@ufrpe.br}
}

\date{Received: August 11, 2017 / Accepted: March 10, 2018}

\maketitle

\begin{abstract}
The consideration of the so-called rotation minimizing frames allows for a simple and elegant characterization of plane and spherical curves in Euclidean space via a linear equation relating the coefficients that dictate the frame motion. In this work, we extend these investigations to characterize curves that lie on a geodesic sphere or totally geodesic hypersurface in a Riemannian manifold of constant curvature. Using that geodesic spherical curves are normal curves, i.e., they are the image of an Euclidean spherical curve under the exponential map, we are able to characterize geodesic spherical curves in hyperbolic spaces and spheres through a non-homogeneous linear equation. Finally, we also show that curves on totally geodesic hypersurfaces, which play the role of hyperplanes in Riemannian geometry, should be characterized by a homogeneous linear equation. In short, our results give interesting and significant similarities between hyperbolic, spherical, and Euclidean geometries.

\end{abstract}

\section{Introduction}
\label{intro}

The geometry of spheres is certainly one of the most important topic of investigation in differential geometry; the search for necessary and/or sufficient conditions for a submanifold be a sphere being one of its major pursuit. A related and interesting problem then is: \emph{how can we characterize those curves $\alpha:I\to \mathbb{R}^{m+1}$ that belong to the surface of a (hyper)sphere?} In $\mathbb{R}^3$, after equipping a curve with its Frenet frame $\{\mathbf{t},\mathbf{n},\mathbf{b}\}$, it is possible to prove that  spherical curves are characterized by the equation
$\kappa/\tau-\mathrm{d}/\mathrm{d}s\left(\kappa'/\tau\kappa^2\right)=0$, where $\kappa$ and $\tau$ are the curvature and torsion, respectively \cite{Kreyszig1991,WongMonatshMath}. Similar relations can be also written in $\mathbb{R}^{m+1}$. On the other hand, by equipping a curve with a {\it rotation minimizing} (RM) {\it frame}, one is able to characterize spherical curves by means of a simple and elegant linear equation involving the coefficients that dictate the frame motion: \emph{a regular curve $\alpha:I\to \mathbb{R}^{m+1}$ is spherical if and only if the normal development curve  $(\kappa_1(s),\dots,\kappa_m(s))$ lies on a line not passing through the origin} \cite{BishopMonthly}. An RM frame $\{\mathbf{t},\mathbf{n}_1,\dots,\mathbf{n}_m\}$ along $\alpha:I\to \mathbb{R}^{m+1}$ is characterized by the equations $\mathbf{t}'(s)=\sum_{i=1}^m\kappa_i(s)\mathbf{n}_i(s)$ and $\mathbf{n}_i'(s)=-\kappa_i(s)\mathbf{t}(s),$
where $s$ is an arc-length parameter. The basic idea here is that $\mathbf{n}_i$ rotates only the necessary amount to remain normal to $\mathbf{t}$: in fact, $\mathbf{n}_i$ is parallel transported along $\alpha$ with respect to the normal connection \cite{Etayo2016}. Due to their minimal twist, RM frames are of importance in applications, such as in computer graphics and visualization \cite{Farouki2008,WangACMTOG2008}, sweep surface modeling \cite{Bloomenthal1991,PottmannIJSM1998,SiltanenCGF1992}, and in differential geometry as well \cite{BishopMonthly,daSilvaArXiv2017,daSilvaJG2017,EtayoTJM2017}, just to name a few.

The goal of this work is to extend these investigations for curves on geodesic spheres in $\mathbb{S}^{m+1}(r)$ and $\mathbb{H}^{m+1}(r)$, the $(m+1)$-dimensional sphere and hyperbolic space of radius $r$, respectively. For spherical curves in $\mathbb{R}^{m+1}$, an important observation is that, up to a translation, their position vectors lie on the normal plane to the curve: $\langle\alpha-p,\alpha-p\rangle=R^2\Leftrightarrow\langle\mathbf{t},\alpha-p\rangle=0$ (we shall call $\alpha$ a \emph{normal curve}). This makes sense due to the double nature of $\mathbb{R}^{m+1}$ as both a manifold and as a tangent space. In fact, this problem has to do with the more general quest of studying curves that lie on a given (moving) plane generated by two chosen vectors of a moving trihedron, e.g., one would define osculating, normal or rectifying curves as those curves whose position vector, up to a translation, lies on their osculating, normal or rectifying planes, respectively \cite{ChenMonthly2003,ChenAJMS2017}: osculating curves are the plane curves (if we substitute the principal normal by an RM vector field, we still have a characterization for plane curves \cite{daSilvaArXiv2017}) and rectifying curves are precisely geodesics on a cone \cite{ChenManuscript,ChenAJMS2017}. This equivalence is no longer valid in other geometries. Nonetheless, it is still possible to extend the concept of normal curves to non-Euclidean settings, such as in affine geometry \cite{KreyszigPAMS1975} and also in $\mathbb{S}^{m+1}(r)$ and $\mathbb{H}^{m+1}(r)$, as we will made clear in this work. Indeed, to extend these notions to a Riemannian setting one should replace the line segment $\alpha(s)-p$ by a geodesic connecting $p$ to $\alpha(s)$, as pointed out by Lucas and Ortega--Yag\"ues in the study of rectifying curves \cite{LucasJMAA2015,LucasMJM2016}: they proved that rectifying curves in the 3d sphere and hyperbolic space are geodesics on a conical surface, in analogy with what happens in the Euclidean case.

Here, we show, as a consequence of the Gauss lemma for the exponential map in a Riemannian manifold $M^{m+1}$, that on a sufficiently small neighborhood of $p\in M^{m+1}$ a curve $\alpha:I\to M^{m+1}$ is normal (with center $p$) if and only if it lies on a geodesic sphere (with center $p$) in $M^{m+1}$. Using this equivalence in the $(m+1)$-dimensional sphere $\mathbb{S}^{m+1}(r)$ and hyperbolic space $\mathbb{H}^{m+1}(r)$, we are able to characterize those curves that lie on the hypersurface of a geodesic sphere in terms of an RM frame. The main result is

{\bf Theorem:} Let $\alpha$ be a regular curve in $\mathbb{S}^{m+1}(r)$ or $\mathbb{H}^{m+1}(r)$. Then, $\alpha$ lies on a geodesic sphere if and only if
\begin{equation}
\left\{\begin{array}{lcc}
\displaystyle\sum_{i=1}^m a_i\,\mathbf{\kappa}_i+\displaystyle\frac{1}{r}\,\cot\left(\frac{z_0}{r}\right) &=0\,, &\mbox{ if }\,\alpha\subseteq\mathbb{S}^{m+1}(r)\\[8pt]
\displaystyle\sum_{i=1}^m a_i\,\mathbf{\kappa}_i+\displaystyle\frac{1}{r}\,\coth\left(\frac{z_0}{r}\right) &=0\,, &\mbox{ if }\,\alpha\subseteq\mathbb{H}^{m+1}(r)\\
\end{array}
\right.,
\end{equation}
for some constants \textcolor{red}{$z_0$} (the radius of the geodesic sphere) and $a_i$, $1\leq i\leq m$.
\newline

For completeness, we also discuss in this work the characterization of geodesic spherical curves in terms of a Frenet frame (Theorem \ref{thr::FrenetChar3DSphCurv}) and show that the characterization of (geodesic) spherical curves is the same as in Euclidean space. Finally, the relation between totally geodesic hypersurfaces, which play the role of hyperplanes in Riemannian geometry, and curves with a normal development   $(\kappa_1,\dots,\kappa_m)$ lying on a line passing through the origin is more delicate, since in general, a manifold has no totally geodesic hypersurfaces up to the trivial ones \cite{MurphyArXiv2017,NikolauevskyIJM2015,Tsukada1996}. Nonetheless, in this work, we are able to show that if a Riemannian manifold contains totally geodesic hypersurfaces, then any curve on a totally geodesic hypersurface is associated with a normal development that lies on a line passing through the origin (Theorem \ref{thr::NormDevelopRiemPlaneCurv}). We show in addition that a curve in $\mathbb{S}^{m+1}(r)$ and $\mathbb{H}^{m+1}(r)$ lies on a totally geodesic hypersurface if and only if its normal development is a line passing through the origin (Theorem \ref{thr::NormDevelopS3andH3PlaneCurv}).

The remaining of this work is organized as follows. In Sect. 2 we review the concept of RM frames, introduce some background material for the geometry of $\mathbb{S}^{m+1}(r)$ and $\mathbb{H}^{m+1}(r)$, and present the concept of normal curves in Riemannian geometry. In Sect. 3 we then characterize geodesic spherical curves via RM and Frenet frames in a constant curvature ambient space, and in Sect. 4, we turn our attention to curves on totally geodesic hypersurfaces. Finally, in Sect. 5, we present our concluding remarks. 

\section{Preliminaries}
 
Let us denote by $\mathbb{E}^{m+1}$ the $(m+1)$-dimensional Euclidean space, i.e., $\mathbb{R}^{m+1}$ equipped with the standard Euclidean metric $\langle\cdot,\cdot\rangle_e$. Given a regular curve $\alpha:I\rightarrow\mathbb{E}^{m+1}$ parametrized by arc-length $s$, i.e., $\langle\mathbf{t},\mathbf{t}\rangle_e=1$, where $\mathbf{t}(s)=\alpha'(s)$, the usual way to introduce a moving frame along it is by means of the Frenet frame $\{\mathbf{e}_0=\mathbf{t},\mathbf{e}_1,\dots,\mathbf{e}_m\}$ \cite{Kreyszig1991,Kuhnel2010}. However, we can also consider any other adapted orthonormal moving frame along $\alpha(s)$: the equation of motion of such a moving frame is then given by a skew-symmetric matrix. Of particular importance are the so-called \emph{Rotation Minimizing (RM) Frames} \cite{BishopMonthly,Etayo2016}: we say that $\{\mathbf{t},\mathbf{n}_1,\dots,\mathbf{n}_m\}$ is an RM frame if $\mathbf{t}$ and $\mathbf{n}_i'$ are parallel. The basic idea here is that $\mathbf{n}_i$ rotates only the necessary amount to remain normal to the tangent $\mathbf{t}$ (so, justifying the terminology). The equation of motion of an RM moving frame is 
\begin{equation}
\frac{{\rm d}}{{\rm d}s}\left(
\begin{array}{c}
\mathbf{t}\\
\mathbf{n}_1\\
\vdots\\
\mathbf{n}_m\\
\end{array}
\right)=\left(
\begin{array}{cccc}
0 & \kappa_{1} & \cdots & \kappa_{m}\\
-\kappa_{1} & 0 & \cdots & 0\\
\vdots & \vdots & \ddots & \vdots \\
-\kappa_{m} & 0 & \cdots & 0\\ 
\end{array}
\right)\left(
\begin{array}{c}
\mathbf{t}\\
\mathbf{n}_1\\
\vdots \\
\mathbf{n}_m\\
\end{array}
\right).\label{eq::BishopEqs}
\end{equation}
\begin{remark}
In $\mathbb{E}^3$, if we write $\mathbf{n}_1=\cos\theta\,\mathbf{n}-\sin\theta\,\mathbf{b}$ and $\mathbf{n}_2=\sin\theta\,\mathbf{n}+\cos\theta\,\mathbf{b}$ for some function $\theta(s)$, the coefficients $\kappa_1,\,\kappa_2$ relate with the curvature function $\kappa$ and torsion $\tau$ according to \cite{BishopMonthly,GuggenheimerCAGD1989}
\begin{equation}
\left\{
\begin{array}{c}
\kappa_1(s) = \kappa(s)\,\cos\theta(s)\\[4pt]
\kappa_2(s) = \kappa(s)\,\sin\theta(s)\\[4pt]
\theta'(s)= \tau(s)\\
\end{array}
\right.\,.\label{eq::RelBetweenRMandFrenet}
\end{equation}
There is a similar relation for curves in $\mathbb{E}^4$ in terms of Euler angles \cite{GokcelikCJMS2014}: we should rotate $\{\mathbf{e}_i\}_{i=1}^3$ to obtain $\{\mathbf{n}_i\}_{i=1}^3$.
\end{remark}

In general, RM frames are not uniquely defined, since any rotation of $\mathbf{n}_i$ on the normal hyperplane still gives an RM field, i.e., there is an ambiguity associated with the action of $SO(m)$ (e.g., in $\mathbb{E}^3$ the angle $\theta$ is only well defined up to an additive constant). Nonetheless, the prescription of curvatures $\kappa_1,\dots,\kappa_m$ uniquely determines a curve up to rigid motions of $\mathbb{E}^{m+1}$ \cite{BishopMonthly,Etayo2016}. In addition, a remarkable advantage of using RM frames is that they allow for a simple characterization of spherical and plane curves:
\begin{theorem}[\cite{BishopMonthly}]
A regular $C^2$ curve $\alpha:I\to \mathbb{E}^{m+1}$ lies on a sphere of radius $r$ if and only if its \emph{normal development}, i.e., the curve $(\kappa_1(s),\dots,\kappa_m(s))$, lies on a line \emph{not passing through the origin}. In addition, $\alpha$ is a plane curve if and only if the normal development lies on a line passing through the origin. \label{thr::charSphCurves}
\end{theorem}

It is also possible to characterize spherical curves through a Frenet frame approach
\begin{theorem}[\cite{Kreyszig1991,Kuhnel2010}]
Let $\alpha:I\to \mathbb{E}^3$ be a $C^4$ regular curve with a non-zero torsion. It lies on a sphere of radius $r$ if and only if
\begin{equation}
\tau(s)\rho(s)+\frac{\mathrm{d}}{\mathrm{d}s}\left(\frac{\rho'(s)}{\tau(s)}\right)=0.
\end{equation}
where $\rho=1/\kappa$ is the radius of curvature.
\end{theorem}
\begin{remark}
It is possible to arrive at a similar characterization for $C^{m+2}$ spherical curves in $\mathbb{E}^{m+1}$, e.g., for spherical curves in $\mathbb{E}^4$ and $\mathbb{E}^5$ with non-zero curvature and torsions, we have
\begin{equation}
\left\{
\begin{array}{c}
\frac{\mathrm{d}}{\mathrm{d}s}\left\{\frac{1}{\tau_2}\frac{\mathrm{d}}{\mathrm{d}s}\left[\frac{1}{\tau_1}\frac{\mathrm{d}}{\mathrm{d}s}\left(\frac{1}{\kappa}\right)\right]+\frac{\tau_1}{\kappa}\right\}+\frac{\tau_2}{\tau_1}\frac{\mathrm{d}}{\mathrm{d}s}\left(\frac{1}{\kappa}\right)=0\\[10pt]
\frac{\mathrm{d}}{\mathrm{d}s}\left\{\frac{1}{\tau_3}\frac{\mathrm{d}}{\mathrm{d}s}\left[\frac{1}{\tau_2}\frac{\mathrm{d}}{\mathrm{d}s}\left[\frac{1}{\tau_1}\frac{\mathrm{d}}{\mathrm{d}s}\frac{1}{\kappa}\right]\right]+\frac{\tau_2}{\tau_1\tau_3}\frac{\mathrm{d}}{\mathrm{d}s}\frac{1}{\kappa}+\frac{1}{\tau_3}\frac{\mathrm{d}}{\mathrm{d}s}\frac{\tau_1}{\kappa\tau_2}\right\}+\frac{\tau_3}{\tau_2}\frac{\mathrm{d}}{\mathrm{d}s}\left[\frac{1}{\tau_1}\frac{\mathrm{d}}{\mathrm{d}s}\frac{1}{\kappa}\right]+\frac{\tau_1\tau_3}{\kappa\tau_2}=0
\end{array}
\right.
\end{equation}
where $\kappa,\tau_1,\dots,\tau_{m-1}$ are the curvature and torsions associated with the Frenet frame $\{\mathbf{e}_i\}_{i=0}^m$: $\mathbf{e}_0'=\tau_0\mathbf{e}_1$ and $\mathbf{e}_i'=-\tau_i\mathbf{e}_{i-1}+\tau_{i+1}\mathbf{e}_{i+1}$ for $1\leq i \leq m$, where $\tau_0=\kappa$ and $\tau_{m+1}=0$ (see Theorem \ref{thr::FrenetChar3DSphCurv}, and comments following it, to have an idea of how devise a proof for the above formulas). Needless to say, the approach via RM frames is simpler, it only demands a $C^2$ condition, and no additional conditions on the torsions and curvature are required.
\end{remark}

\subsection{Rotation minimizing frames and normal curves in Riemannian geometry}

It is also possible to introduce Frenet frames in Riemannian manifolds \cite{GutkinJGP2011,Spivak1979v4}, see also \cite{BolcskeiBAG2007,GutkinJGP2011,LucasJMAA2015,LucasMJM2016,SzilagyiSUZ2003}. Analogously, one can also define RM frames \cite{Etayo2016,EtayoTJM2017}. To introduce such concepts, one should take covariant derivatives in the direction of the unit tangent instead of the ordinary one. More precisely, let $M^{m+1}$ be a Riemannian manifold with Levi--Civita connection $\nabla$ and metric $\langle\cdot,\cdot\rangle$ \cite{doCarmo1992}. We say that $\mathbf{x}\in\mathfrak{X}(M)$ is an RM vector field along a regular curve $\alpha:I\to M^{m+1}$ if $\nabla_{\mathbf{t}}\,\mathbf{x}=\lambda\,\mathbf{t}$, where $\mathfrak{X}(M)$ is the module of tangent vector fields, $\mathbf{t}(s)=\alpha'(s)$ is the unit tangent, and $s$ an arc-length parameter \cite{Etayo2016}.

To build a Frenet frame in $M$, the curvature function and principal normal (if $\kappa\not=0$) are defined as usual, that is
\begin{equation}
\kappa =\Vert\nabla_{\mathbf{t}}\,\mathbf{t}\,\Vert\,\,\mbox{ and }\,\,\mathbf{n} = \frac{1}{\kappa}\nabla_{\mathbf{t}}\,\mathbf{t}\,,
\end{equation}
respectively. The binormal vector $\mathbf{b}$ is chosen in a way that $\{\mathbf{t},\mathbf{n},\mathbf{b}\}$ is a positively oriented orthonormal frame along $T_{\alpha(s)}M$. The torsion is given by
\begin{equation}
\tau = -\langle\nabla_{\mathbf{t}}\,\mathbf{b},\mathbf{n}\,\rangle\,,
\end{equation}
and the Frenet equations can be written as
\begin{equation}
\nabla_{\mathbf{t}}\left(
\begin{array}{c}
\mathbf{t}\\
\mathbf{n}\\
\mathbf{b}\\
\end{array}
\right)=\left(
\begin{array}{ccc}
0 & \kappa & 0\\
-\kappa & 0 & \tau\\
 0 & -\tau & 0\\ 
\end{array}
\right)\left(
\begin{array}{c}
\mathbf{t}\\
\mathbf{n}\\
\mathbf{b}\\
\end{array}
\right).\label{eq::RiemFrenetEqs}
\end{equation}

In this work, we will be primarily interested in the $(m+1)$-dimensional sphere $\mathbb{S}^{m+1}(r)$ and in the hyperbolic space $\mathbb{H}^{m+1}(r)$. We will, respectively, use them modeled as submanifolds of $\mathbb{E}^{m+2}$ and $\mathbb{E}_1^{m+2}$:
\begin{equation}
\mathbb{S}^{m+1}(r) = \{q\in\mathbb{R}^{m+2}\,:\,\langle q,q\rangle_e=r^2\}
\end{equation}
and
\begin{equation}
\mathbb{H}^{m+1}(r) = \{q\in\mathbb{R}^{m+2}\,:\,\langle q,q\rangle_1=-r^2,\,x_1>0\},
\end{equation}
equipped with the induced metric denoted by $\langle\cdot,\cdot\rangle$ (the context will make clear if we are using $\langle\cdot,\cdot\rangle_e$ or $\langle\cdot,\cdot\rangle_1$). Here, $\mathbb{E}_1^{m+2}$ denotes the Lorentz space equipped with the index 1 metric $\langle \mathbf{x},\mathbf{y}\rangle_1=-x_1y_1+\sum_{i=2}^{m+2}x_iy_i$. 

Denoting by $\nabla$ and $\nabla^0$ the Levi--Civita connections on $\mathbb{S}^{m+1}(r)$ (or $\mathbb{H}^{m+1}(r)$) and $\mathbb{E}^{m+2}$ (or $\mathbb{E}^{m+2}_1$, respectively), they are related by the Gauss formula as follows:
\begin{equation}
\nabla^0_{\mathbf{x}}\,\mathbf{y}=\nabla_{\mathbf{x}}\,\mathbf{y}\mp\frac{1}{r^2}\langle\mathbf{x},\mathbf{y}\rangle \,q\,,\label{eq::CovDerAndUsualDer}
\end{equation}
where $q$ denotes the position vector, i.e., the canonical immersion $q:\mathbb{S}^{m+1}(r)\to\mathbb{E}^{m+2}$ for the minus sign and $q:\mathbb{H}^{m+1}(r)\to\mathbb{E}^{m+2}_1$ for the plus sign.  

\begin{remark}
The models above do not represent the unique choices. Another common way of looking at the spherical geometry is the intrinsic model based on stereographic projection \cite{doCarmo1992,Spivak1979v4}. On the other hand, besides the hyperboloid model above, other common models for the hyperbolic space are the Poincar\'e ball and half-plane models \cite{BenedettiPetronio,ReynoldsMonthly1993,Spivak1979v4}. In any case, the important fact is that these models are all isometric. Thus, intrinsically speaking, they are all the same, and the choice between them being a matter of convenience.
\end{remark}

The concept of normal curves will be of fundamental importance in our work. In Euclidean space we say that $\alpha$ is a \emph{normal curve} if
\begin{equation}
\alpha(s)-p \in \mbox{span}\{\mathbf{t}(s)\}^{\perp},
\end{equation}
where $p$ is a fixed point (the center of the normal curve). We can straightforwardly prove that normal curves in $\mathbb{E}^{m+1}$ are precisely the spherical ones (in this case, $p$ is the center of the respective sphere): $\langle \alpha-p,\mathbf{t}\rangle=0\Leftrightarrow \langle \alpha-p,\alpha-p\rangle=$ constant. This definition makes sense due to the double nature of $\mathbb{E}^{m+1}$ as both a manifold and a tangent space. To extend it to a Riemannian manifold $M^{m+1}$, we should replace $\alpha-p$ by a geodesic connecting $p$ to a point $\alpha(s)$ on the curve, as done in \cite{LucasJMAA2015,LucasMJM2016} for the study of rectifying curves:
\begin{definition}
A regular curve $\alpha:I\to \mathbb{S}^{m+1}(r)$ or $\alpha:I\to \mathbb{H}^{m+1}(r)$ is a \emph{normal curve} with center $p$ if the geodesic $\beta_s$ connecting $p$ to $\alpha(s)$ is orthogonal to $\alpha$, i.e., $\langle\mathbf{t}_{\alpha},\mathbf{t}_{\beta}\rangle (s)=0$ for all $s\in I$. In $\mathbb{S}^{m+1}(r)$ we additionally assume that $\alpha$ does not contain the antipodal $-p$ of $p$: $-p\not\in\mbox{Im}(\alpha)$.
\end{definition}

The above definition is also valid in a generic Riemannian manifold $M^{m+1}$ once we restrict ourselves to work on a sufficiently small neighborhood of $p$ (out of the injectivity radius the geodesic $\beta_s$ may fail to be unique). The equivalence between spherical and normal curves can be extended to a Riemannian manifold by applying the Gauss lemma for the exponential map \cite{doCarmo1992}:

\begin{proposition}
On a sufficiently small neighborhood of $p\in M^{m+1}$, a curve $\alpha:I\to M^{m+1}$ is normal (with center $p$) if and only if it lies on a geodesic sphere (with center $p$). In other words, a normal curve is the image of an Euclidean spherical curve under the exponential map.\label{prop::EquivNormalAndSphCurves}
\end{proposition}

Finally, given $p\in\mathbb{S}^{m+1}(r)$, $\mathbf{v}\in \mathbb{S}^m(1)\subset T_p\mathbb{S}^{m+1}(r)$ or $p\in\mathbb{H}^{m+1}(r)$, $\mathbf{v}\in\mathbb{S}^m(1)\subset T_p\mathbb{H}^{m+1}(r)$, the exponential map is
\begin{equation}
\mbox{exp}_p(u\mathbf{v})=\cos\left(\frac{u}{r}\right)p+r\sin\left(\frac{u}{r}\right)\mathbf{v}\,
\end{equation}
or
\begin{equation}
\mbox{exp}_p(u\mathbf{v})=\cosh\left(\frac{u}{r}\right)p+r\sinh\left(\frac{u}{r}\right)\mathbf{v}\,,
\end{equation}
respectively. Observe that the geodesics $\beta(u)=\exp_p(u\mathbf{v})$ above are defined for any value $u\in\mathbb{R}$ and then the equivalence in Proposition \ref{prop::EquivNormalAndSphCurves} is valid globally.

\section{Spherical curves in $\mathbb{S}^{m+1}(r)$ and $\mathbb{H}^{m+1}(r)$}

\begin{figure*}[tbp]
\centering
    {\includegraphics[width=0.33\linewidth]{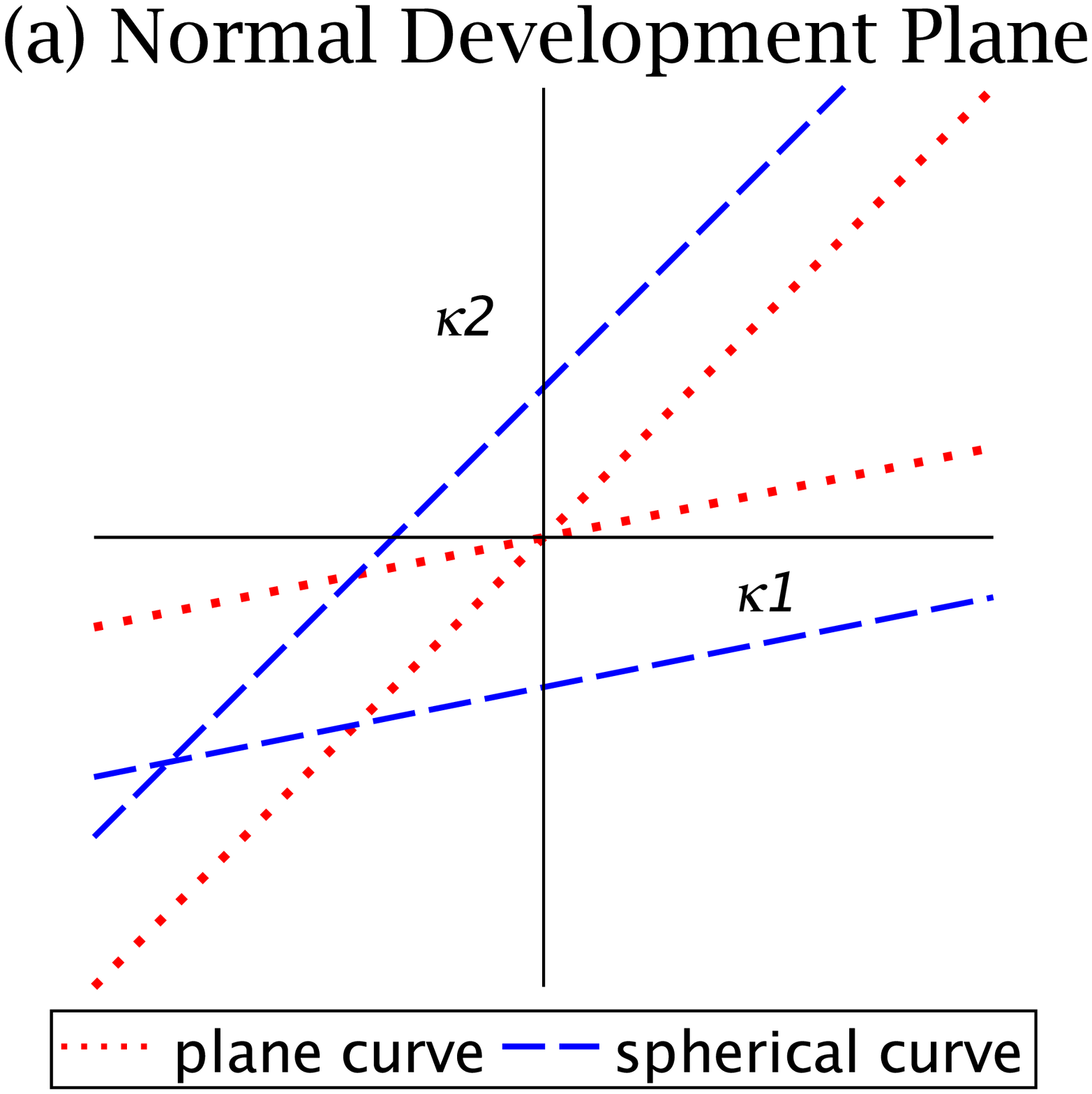}}
    {\includegraphics[width=0.32\linewidth]{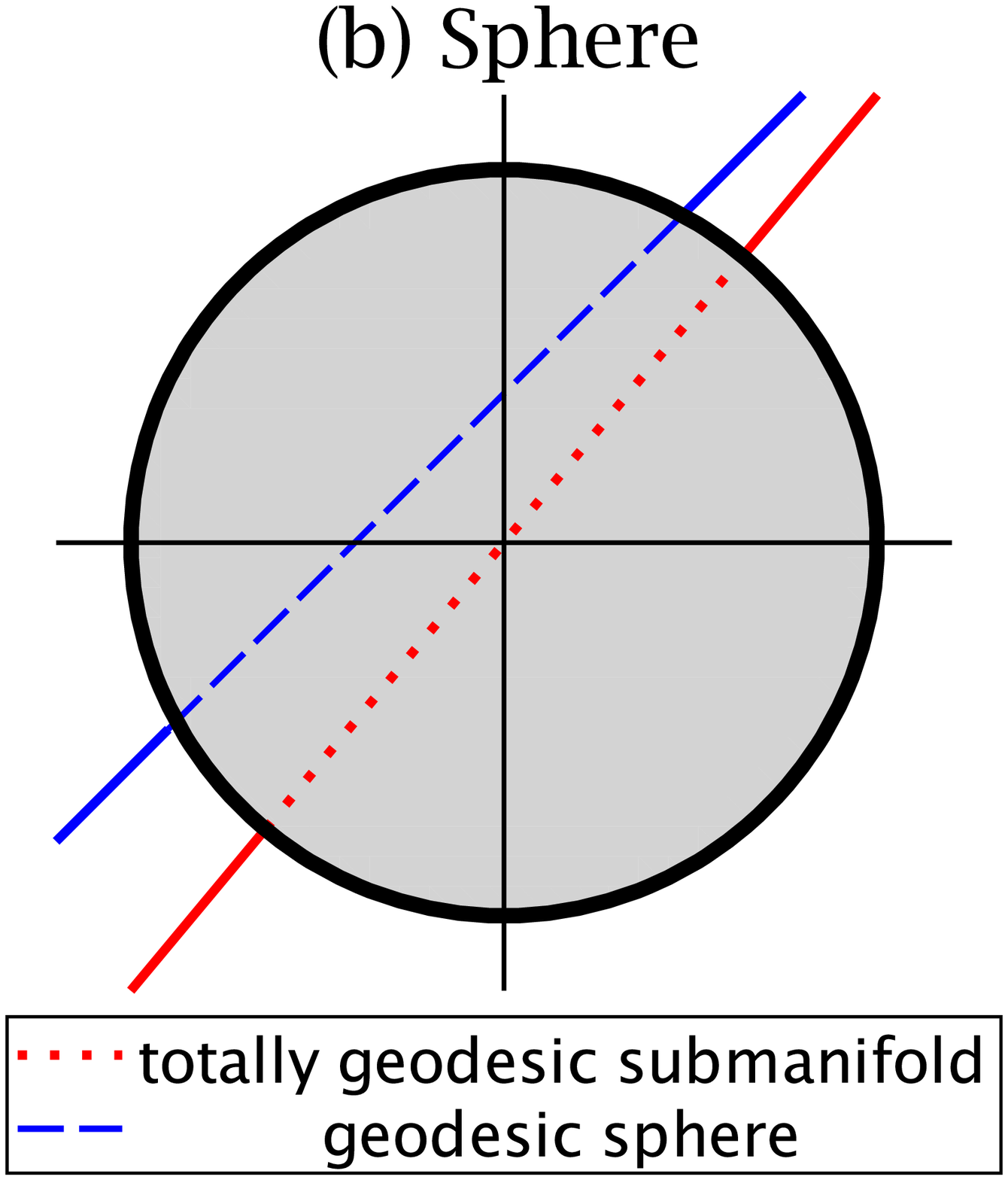}}
    {\includegraphics[width=0.32\linewidth]{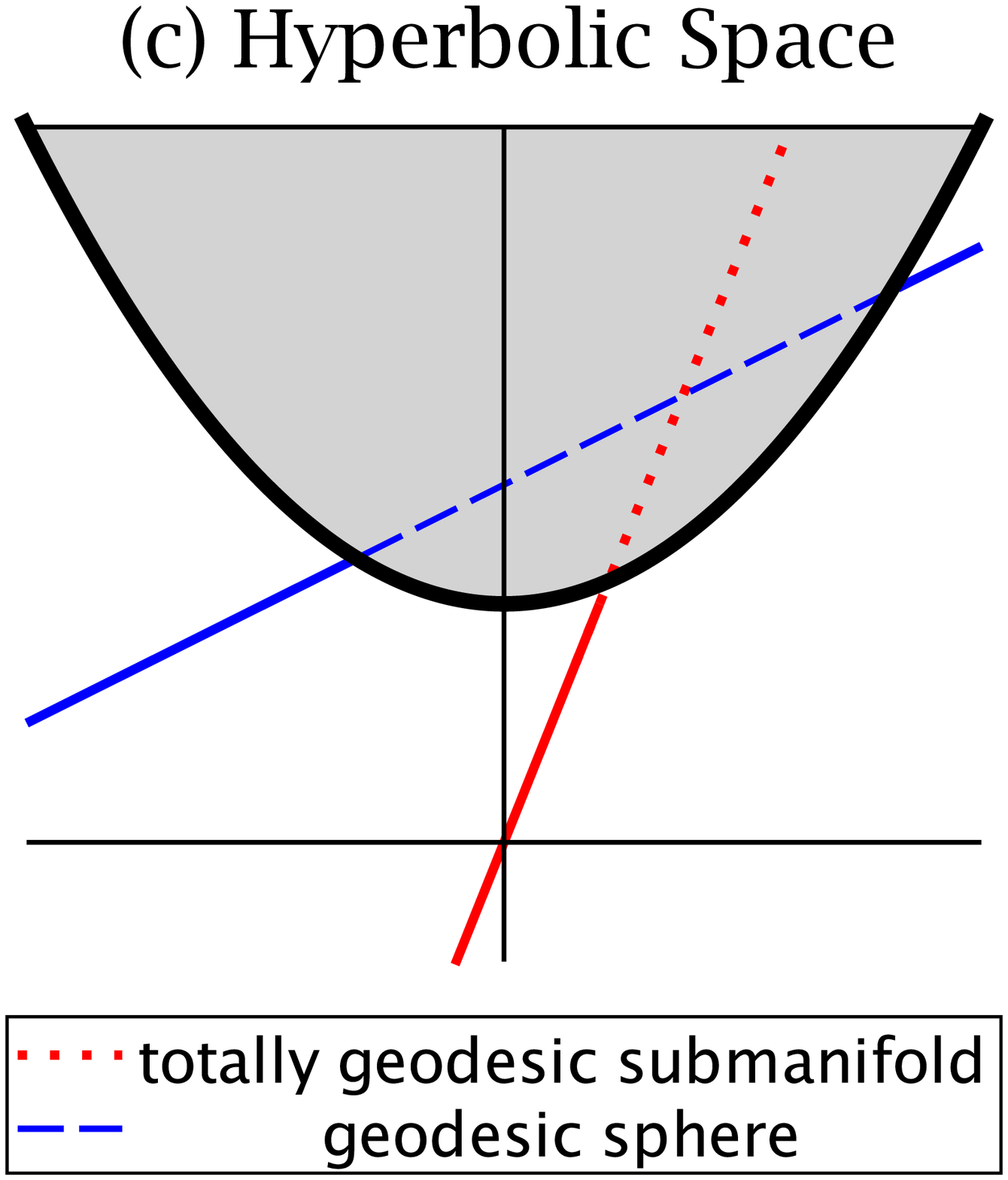}}
          \caption{The geometry of the normal development, geodesic spheres, and totally geodesics submanifolds in $\mathbb{S}^{m+1}(r)$ and $\mathbb{H}^{m+1}(r)$: \textbf{(a)} ($m=2$ in the figure) lines not passing through the origin (dashed blue line) represent geodesic spherical curves (Theorem \ref{thr::RMCharSphGeodCurves}) and lines through the origin (dotted red line) represent plane curves, i.e., curves on totally geodesic hypersurfaces (Theorem \ref{thr::NormDevelopS3andH3PlaneCurv}); \textbf{(b)} and \textbf{(c)} lines passing through the origin (dotted red line) represent hyperplanes passing through the origin and, when intersected with $\mathbb{S}^{m+1}(r)$  or $\mathbb{H}^{m+1}(r)$, give rise to totally geodesic hypersurfaces, while lines not passing through the origin (dashed blue line) represent hyperplanes not passing through the origin and when intersected with $\mathbb{S}^{m+1}(r)$ or $\mathbb{H}^{m+1}(r)$, give rise to geodesic spheres (in the hyperboloid model, an intersection with hyperplanes forming smaller angles with the hyperboloid axis give rise to equidistant surfaces and horospheres) \cite{Spivak1979v4}.}
          \label{fig::DiagramSphPlaneCurv}
\end{figure*}

As previously said, by equipping a curve in $\mathbb{E}^{m+1}$ with an RM frame, it is possible to characterize spherical curves by means of a linear relation involving the coefficients which dictate the frame motion. We now extend these results for curves on geodesic spheres of $\mathbb{S}^{m+1}(r)$ and $\mathbb{H}^{m+1}(r)$ (see Fig. \ref{fig::DiagramSphPlaneCurv}).
\begin{theorem}
Let $\alpha$ be a regular $C^2$ curve in $\mathbb{S}^{m+1}(r)$ or $\mathbb{H}^{m+1}(r)$. Then, $\alpha$ lies on a geodesic sphere if and only if
\begin{equation}
\left\{\begin{array}{lcc}
\displaystyle\sum_{i=1}^m a_i\,\mathbf{\kappa}_i+\displaystyle\frac{1}{r}\,\cot\left(\frac{z_0}{r}\right) &=0\,, &\mbox{ if }\,\alpha\subseteq\mathbb{S}^{m+1}(r)\\[5pt]
\displaystyle\sum_{i=1}^m a_i\,\mathbf{\kappa}_i+\displaystyle\frac{1}{r}\,\coth\left(\frac{z_0}{r}\right) &=0\,, &\mbox{ if }\,\alpha\subseteq\mathbb{H}^{m+1}(r)\\
\end{array}
\right.,
\end{equation}
for some constants $z_0$ (the radius of the geodesic sphere\footnote{For $\mathbb{S}^{m+1}(r)$ we may impose $z_0<\pi r/2$, which guarantees that the center of the geodesic sphere is well defined: if $z_0=\pi r/2$, both $p$ and its antipodal $-p$ are equidistant from the geodesic sphere.}) and $a_i$, $1\leq i\leq m$. \label{thr::RMCharSphGeodCurves}
\end{theorem}
{\it Proof. } We will do the proof for $\alpha\subseteq\mathbb{S}^{m+1}(r)$ only, the case for $\mathbb{H}^{m+1}(r)$ being analogous (one just needs to use the hyperbolic versions of the trigonometric functions). 

If $\alpha:I\to\mathbb{S}^{m+1}(r)$ is a normal curve parametrized by arc-length $s$, then we may write
\begin{equation}
\alpha(s) =  \exp_p(z_0\,\mathbf{v}(c_0\,s)),
\end{equation}
where $z_0$ and $c_0=[r\sin(z_0/r)]^{-1}$ are constants and $\mathbf{v}:I\to\mathbb{S}^{m}(1)\subseteq T_p\mathbb{S}^{m+1}(r)$ is a unit speed curve. In our model of $\mathbb{S}^{m+1}(r)$ as submanifold of $\mathbb{E}^{m+2}$, a tangent vector $\mathbf{v}$ at $p$ satisfies $\langle p,\mathbf{v}\rangle_e=0$. If $\{\mathbf{t}_{\alpha},\mathbf{n}_1,\dots,\mathbf{n}_m\}$ is an RM frame along $\alpha$, the unit tangent of $\alpha$ can be written as
\begin{equation}
\mathbf{t}_{\alpha}(s) = \mathbf{v}'(c_0s).\label{eqTalphaInTermsOfV}
\end{equation}

On the other hand, the unit speed geodesic $\beta_s$ connecting $p$ to a point $\alpha(s)$ is
\begin{equation}
\beta_s(u) = \cos\left(\frac{u}{r}\right)\,p+r\sin\left(\frac{u}{r}\right)\,\mathbf{v}(c_0\,s)\,\Rightarrow \beta_s(z_0)=\alpha(s).
\end{equation}

The normality condition $\langle \mathbf{t}_{\alpha},\mathbf{t}_{\beta}\rangle=0$ implies
\begin{equation}
``\mathbf{t}_{\beta_s}\mbox{ at }\alpha(s)" = \mathbf{t}_{\beta_s}(z_0)=\displaystyle\sum_{i=1}^ma_i(s)\mathbf{n}_i(s)\,.
\end{equation}

The derivative of the coefficients $a_i=\langle \mathbf{t}_{\beta_s},\mathbf{n}_i\rangle$ gives
\begin{equation}
a'_i = \langle\nabla_{\mathbf{t}_{\alpha}}\mathbf{t}_{\beta_s},\mathbf{n}_i\rangle+\langle \mathbf{t}_{\beta_s},\nabla_{\mathbf{t}_{\alpha}}\mathbf{n}_i\rangle= \langle \nabla^0_{\mathbf{t}_{\alpha}}\mathbf{t}_{\beta_s},\mathbf{n}_i\rangle,\,\label{eq::DerivativeOfai}
\end{equation}
where the last equality is a consequence of the fact that $\mathbf{n}_i$ is RM and also that $\nabla_{\mathbf{x}}\,\mathbf{y}=\nabla_{\mathbf{x}}^0\,\mathbf{y}$ for two orthogonal vectors $\mathbf{x},\,\mathbf{y}$, see Eq. (\ref{eq::CovDerAndUsualDer}). Now, using that $\mathbf{t}_{\beta}$ along $\alpha$ can be also written as
\begin{equation}
\mathbf{t}_{\beta_s}(z_0) = -\frac{1}{r}\sin\left(\frac{z_0}{r}\right)\,p+\cos\left(\frac{z_0}{r}\right)\,\mathbf{v}(c_0\,s),\label{eq::TbetaAlongAlpha}
\end{equation}
we have
\begin{equation}
\nabla_{\mathbf{t}_{\alpha}}^0\mathbf{t}_{\beta} = \frac{1}{r}\frac{\cos(z_0/r)}{\sin(z_0/r)}\,\mathbf{v}'(c_0\,s)=\frac{\cot(z_0/r)}{r}\,\mathbf{t}_{\alpha}\,.\label{eq::NablaTbAlongTaIsTa}
\end{equation}
Inserting the expression above in Eq. (\ref{eq::DerivativeOfai}) shows that $a_i'=0$, and therefore, the coefficients $a_i$, $1\leq i\leq m$, are all constants.

Finally, taking the derivative of $\langle \mathbf{t}_{\beta},\mathbf{t}_{\alpha}\rangle=0$ along $\alpha$ gives
\begin{eqnarray}
0 & = & \langle\nabla_{\mathbf{t}_{\alpha}}\mathbf{t}_{\beta},\mathbf{t}_{\alpha}\rangle+\langle \mathbf{t}_{\beta},\nabla_{\mathbf{t}_{\alpha}}\mathbf{t}_{\alpha}\rangle\nonumber\\
& = & \left\langle\frac{\cot(z_0/r)}{r}\,\mathbf{t}_{\alpha},\mathbf{t}_{\alpha}\right\rangle+\left\langle \sum_{i=1}^ma_i\mathbf{n}_i,\sum_{j=1}^m\kappa_j\mathbf{n}_j\right\rangle\nonumber\\
& = & \frac{1}{r}\cot\left(\frac{z_0}{r}\right)+\sum_{i=1}^ma_i\kappa_i.
\end{eqnarray}

Conversely, suppose that $\alpha$ is a regular curve and that it satisfies $\sum_{i}a_i\kappa_i+\cot(z_0r^{-1})/r=0$. The proof is based on the following observation: for a spherical curve, if we invert the direction of the motion of $\beta_s$ we have a geodesic connecting $\alpha(s)$ to $p$, whose initial velocity vector according to Eq. (\ref{eq::TbetaAlongAlpha}) should be $-\mathbf{t}_{\beta}$. Now, let us define
\begin{equation}
\mathbf{w}(s)=-\sum_{i=1}^ma_i\mathbf{n}_i
\end{equation}
and
\begin{equation}
P(s) = \cos\left(\frac{z_0}{r}\right)\alpha(s) - r\sin\left(\frac{z_0}{r}\right)\mathbf{w}(s).
\end{equation}
Taking the derivative of the last equation, we find $P'(s)=0$ and then $P$ is a constant point. Consequently, it means that the geodesics with initial point $\alpha(s)$ and initial velocity $\mathbf{w}(s)$ travel always the same distance to arrive at $P$, i.e., $\alpha$ is a spherical curve.

\qed

Finding RM frames along a curve may be a difficult problem and, in general, one must resort to some kind of numerical method, see e.g. \cite{WangACMTOG2008}. However, for a curve $\alpha$ in $\mathbb{S}^2(r,p)\subseteq\mathbb{R}^3$, computing RM frames is not difficult: $\mathbf{u}=(\alpha(s)-p)/r$ is RM \cite{daSilvaArXiv2017,WangACMTOG2008}. This result can be extended for other ambient spaces by taking into account Eq. (\ref{eq::NablaTbAlongTaIsTa}) in the proof above. Then, we have

\begin{corollary}
For a regular $C^2$ curve $\alpha(s)$ on a geodesic sphere of $\mathbb{S}^{m+1}(r)$, or $\mathbb{H}^{m+1}(r)$, the tangents of the geodesics connecting the center of the geodesic sphere to points on the curve is a rotation minimizing vector field. 
\end{corollary}

The previous theorem was obtained by expressing $\mathbf{t}_{\beta}$ in terms of an RM basis for the normal plane $\mbox{span}\{\mathbf{t}_{\alpha}\}^{\perp}$. If we use the Frenet frame instead, then we can extend a classical characterization result for spherical curves in $\mathbb{R}^3$.

\begin{theorem}
Let $\alpha$ be a regular $C^4$ curve with non-zero torsion in $\mathbb{S}^3(r)$ or $\mathbb{H}^3(r)$. The curve $\alpha$ lies on a geodesic sphere if, and only if
\begin{equation}
\frac{{\rm d}}{{\rm d}s}\left[\frac{1}{\tau}\frac{{\rm d}}{{\rm d}s}\left(\frac{1}{\kappa}\right)\right]+\frac{\tau}{\kappa}=0\,.\label{eq::FrenetCharSphCurves}
\end{equation}
\label{thr::FrenetChar3DSphCurv}
\end{theorem}
{\it Proof. } We will do the proof for $\alpha\subseteq\mathbb{S}^3(r)$ only, the case for $\mathbb{H}^3(r)$ being analogous. 

Let $\alpha$ be a spherical curve and $\{\mathbf{t}_{\alpha},\mathbf{n},\mathbf{b}\}$ its Frenet frame, then there exists a point $p$ such that the geodesic $\beta_s$ connecting $p$ to $\alpha(s)$ satisfies $\langle\mathbf{t}_{\beta},\mathbf{t}_{\alpha}\rangle=0$. Let us write
\begin{equation}
\mathbf{t}_{\beta} = c_1\mathbf{n}+c_2\mathbf{b},
\end{equation}
for some functions $c_1,\,c_2$.

Taking the (covariant) derivative gives
\begin{eqnarray}
\nabla_{\mathbf{t}_{\alpha}}\mathbf{t}_{\beta} & = & c'_1\mathbf{n}+c'_2\mathbf{b}+c_1\nabla_{\mathbf{t}_{\alpha}}\mathbf{n}+c_2\nabla_{\mathbf{t}_{\alpha}}\mathbf{b}\nonumber\\
\frac{1}{r}\cot\left(\frac{z_0}{r}\right)\mathbf{t}_{\alpha} & = & -c_1\kappa\mathbf{t}_{\alpha}+(c'_1-\tau c_2)\mathbf{n}+(c'_2+\tau c_1)\mathbf{b},
\end{eqnarray}
where we used Eqs. (\ref{eq::RiemFrenetEqs}) and (\ref{eq::NablaTbAlongTaIsTa}) to arrive at the second equality above. Now, comparing coefficients leads to 
\begin{equation}
\left\{
\begin{array}{ccc}
-\kappa c_1 & = & \frac{1}{r}\cot\left(\frac{z_0}{r}\right)\\
c'_1 - \tau c_2 & = & 0\\
c'_2 + \tau c_1 & = & 0\\
\end{array}
\right..\label{eq::auxFrenetCharSphCurves}
\end{equation}
 
 From the first and second equations, we find
\begin{equation}
c_1=-\frac{1}{r\kappa}\cot\left(\frac{z_0}{r}\right)\Rightarrow \tau c_2 = c'_1 = -\frac{{\rm d}}{{\rm d}s}\left[\frac{1}{r\kappa}\cot\left(\frac{z_0}{r}\right)\right].
\end{equation}
Now, using the expression above in combination with the 3rd equation of (\ref{eq::auxFrenetCharSphCurves}) furnishes
\begin{equation}
-\tau c_1 = c_2'=-\frac{{\rm d}}{{\rm d}s}\left\{\frac{1}{\tau}\frac{{\rm d}}{{\rm d}s}\left[\frac{1}{r\kappa}\cot\left(\frac{z_0}{r}\right)\right]\right\}\,.
\end{equation}
The desired result follows from the finding above and the 1st equation of (\ref{eq::auxFrenetCharSphCurves}).

Conversely, let $\alpha$ be a regular curve satisfying Eq. (\ref{eq::FrenetCharSphCurves}). As in the proof for the characterization of spherical curve via RM frames, the idea is to find a (fixed) point $P$ and a vector field $\mathbf{w}$ such that all the geodesics emanating from $\alpha$ with initial velocity $\mathbf{w}$ reach $P$ after traveling the same distance. Let us define the following vector field along $\alpha(s)$
\begin{equation}
\mathbf{w}(s)= -\frac{1}{r\kappa(s)}\cot\left(\frac{z_0}{r}\right)\mathbf{n}(s)-\frac{1}{\tau(s)}\frac{\textrm{d}}{\textrm{d}s}\left[\frac{1}{r\kappa(s)}\cot\left(\frac{z_0}{r}\right)\right]\,\mathbf{b}(s),
\end{equation}
which satisfies $\nabla_{\mathbf{t}_{\alpha}}\mathbf{w}=r^{-1}\cot(z_0r^{-1})\,\mathbf{t}_{\alpha}\,$. Now, define
\begin{equation}
P(s) = \cos\left(\frac{z_0}{r}\right)\alpha(s)-r\sin\left(\frac{z_0}{r}\right)\mathbf{w}(s)\,.
\end{equation}
Taking the derivative of $P$ shows that $P'(s)=0$, and therefore, $P$ is constant and will be the center of the geodesic sphere that contains $\alpha$.
\qed

\begin{remark}
One can also equip a curve with a Frenet frame in higher dimensional Riemannian manifolds \cite{Spivak1979v4}, p. 29, and use them to characterize (geodesic) spherical curves. One can follow the same steps as in the previous theorem, i.e., use that a spherical curve must be normal and then investigate the coefficients $c_i$ of $\mathbf{t}_{\beta}$ in terms of the Frenet frame. The expressions, however, are quite cumbersome and we will not attempt to write it here. We just remark that, as happens in 3d, the values of $r$ and of the geodesic sphere radius do not appear in the expression characterizing spherical curves. Note in addition that the curve must be of class $C^{m+2}$, in contrast with the $C^2$ requirement in Theorem \ref{thr::RMCharSphGeodCurves} via RM frames.
\end{remark}

\section{Curves on totally geodesic hypersurface}

The so-called totally geodesic submanifolds in a Riemannian ambient space have the simplest shape and play the role of affine subspaces. Despite their simplicity, in general, Riemannian manifolds do not have non-trivial totally geodesic submanifolds \cite{MurphyArXiv2017,Tsukada1996}. The existence of such submanifolds imposes severe restrictions on the geometry of the ambient manifold \cite{NikolauevskyIJM2015}. Riemannian space forms are examples of manifolds that contain non-trivial totally geodesic submanifolds.

\begin{definition}
A submanifold $N$ of a Riemannian manifold $M$ is a \emph{totally geodesic submanifold} if any geodesic on the submanifold $N$ with the induced Riemannian metric is also a geodesic on $M$ (e.g., one dimensional totally geodesic submanifolds are geodesics).
\end{definition}

In the following, we shall restrict our attention to orientable hypersurfaces. There are many equivalent ways of characterizing a totally geodesic hypersurface. Indeed, all the conditions below are equivalent \cite{Cartan1946}, p. 114,
\begin{enumerate}
\item $N\subset M$ is totally geodesic.
\item the principal curvatures vanish in every point of $N$.
\item the normal field to $N$ remains normal if parallel transported along any curve on $N$.
\item any tangent field to $N$ remains tangent if parallel transported along any curve on $N$.
\end{enumerate}
Note that property 3 essentially says that the normal field of a totally geodesic hypersurface is constant, which is a crucial feature of Euclidean hyperplanes: $\pi$ is a hyperplane if and only if there exist $\mathbf{u}_0$ and $x_0$ constants, such that $\pi=\{x:\langle x-x_0,\mathbf{u}_0\rangle=0\}$. Thus, we may see \emph{hyperplane curves} in a Riemannian manifold as those curves on totally geodesic hypersurfaces.

In Euclidean space, it is known that normal development curves $(\kappa_1,\dots,\kappa_m)$ which are lines passing through the origin characterize hyperplane curves (Theorem \ref{thr::charSphCurves}). Here, we (partially) extend this result to totally geodesic curves on \emph{any} Riemannian manifold.
\begin{theorem}
Let $\alpha:I\to N^m\subset M^{m+1}$ be a regular curve and $\{\mathbf{t},\mathbf{n}_1,\dots,\mathbf{n}_m\}$ a rotation minimizing frame along it. If $\alpha$ lies on a totally geodesic hypersurface $N$, then its normal development curve 
$(\kappa_1,\dots,\kappa_m)$ lies on a line passing through the origin.
\label{thr::NormDevelopRiemPlaneCurv}
\end{theorem}
{\it Proof. } Let $\mathbf{u}$ be a normal vector field on $N$. Since $N$ is totally geodesic, we can use that $\nabla_{\mathbf{t}}\,\mathbf{u}=0$. In addition, we can also write $\mathbf{u}=\sum_{i=1}^ma_i\mathbf{n}_i$ for the normal $\mathbf{u}$ along $\alpha$. The coefficient $a_i=\langle\mathbf{u},\mathbf{n}_i\rangle$ satisfies
\begin{equation}
a'_i = \langle\nabla_{\mathbf{t}}\,\mathbf{u},\mathbf{n}_i\rangle+\langle\mathbf{u},\nabla_{\mathbf{t}}\,\mathbf{n}_i\rangle=0.
\end{equation}
Then, for all $i\in\{1,\dots,m\}$, $a_i$ is a constant. Finally
\begin{equation}
0=\nabla_{\mathbf{t}}\,\mathbf{u}=\sum_{i=1}^ma_i\nabla_{\mathbf{t}}\,\mathbf{n}_i=\sum_{i=1}^m(-a_i\kappa_i\,\mathbf{t})
\end{equation}
and therefore, $\sum a_i\kappa_i=0$ represents the equation of a line passing through the origin.
\qed

Let us now discuss the reciprocal of the theorem above. Given a curve $\alpha:I\to M$ satisfying $\sum_{i=1}^ma_i\,\kappa_i=0$ for some constants $a_1,\dots,a_m$, we may define $\mathbf{u}(s)=\sum a_i\,\mathbf{n}_i$. Then, it follows that
\begin{equation}
\nabla_{\mathbf{t}}\,\mathbf{u}=\sum -a_i\kappa_i\,\mathbf{t}=0.
\end{equation}
Thus, $\mathbf{u}$ is parallel transported along $\alpha$. The problem now is to find a codimension 1 totally geodesic submanifold containing $\alpha$ and whose normal field is equal to $\mathbf{u}$ when restricted to $\alpha$. A candidate to solution is the submanifold given by the following parametrization:
\begin{equation}
X(s_1,\dots,s_m) = \exp_{\alpha(s_1)}\left(\sum_{i=2}^{m}s_i\,\mathbf{u}_i\right),
\end{equation}
where $\{\mathbf{u}_i(s)\}_{i=2}^{m}$ is an orthonormal basis for  $\mbox{span}\{\mathbf{t}(s),\mathbf{u}(s)\}^{\perp}$ for all $s=s_1$. Observe, however, the fact that $X$ is geodesic along $\alpha$ does not implies that it will also be geodesic in all its points. In fact, the existence of non-trivial totally geodesic submanifolds is an exceptional fact. On the other hand, in both $\mathbb{S}^{m+1}(r)$ and $\mathbb{H}^{m+1}(r)$ the situation is easier, since that totally geodesic submanifolds do exist and are precisely the intersection of affine subspaces of $\mathbb{R}^{m+2}$ with $\mathbb{S}^{m+1}(r)$ and $\mathbb{H}^{m+1}(r)$ \cite{Spivak1979v4} (see Fig. \ref{fig::DiagramSphPlaneCurv}). Then, we have
\begin{theorem}
Let $\alpha$ be a regular $C^2$ curve in $\mathbb{S}^{m+1}(r)$, or $\mathbb{H}^{m+1}(r)$, equipped with an RM frame $\{\mathbf{t},\mathbf{n}_1,\dots,\mathbf{n}_m\}$. Then, $\alpha$ is a hyperplane curve, i.e., it lies on a totally geodesic hypersurface, if and only if the normal development $(\kappa_1,\dots,\kappa_m)$ is a line passing through the origin.
\label{thr::NormDevelopS3andH3PlaneCurv}
\end{theorem}
\textit{Proof.} The direction ``hyperplane curve $\Rightarrow$ $\sum_{i=1}^ma_i\kappa_i=0$ ($a_i$ constant)'' is a consequence of the previous theorem. For the reciprocal, define a vector field along $\alpha$ as $\mathbf{u}(s)=\sum_{i=1}^ma_i\mathbf{n}_i(s)=0$. Using that the normal development is a line passing through the origin, we have 
\begin{equation}
\frac{\mathrm{d}\mathbf{u}}{\mathrm{d}s}\equiv\nabla^0_{\mathbf{t}}\,\mathbf{u}=\nabla_{\mathbf{t}}\,\mathbf{u} = \sum_{i=1}^m-a_i\kappa_i\,\mathbf{t} = 0,
\end{equation}
where for the second equality, we used that $\langle\mathbf{t},\mathbf{u}\rangle=0$ in Eq. (\ref{eq::CovDerAndUsualDer}). Therefore, $\mathbf{u}$ is a constant vector in $\mathbb{R}^{m+2}$ and it follows that $\alpha$ is contained in the hyperplane, in $\mathbb{R}^{m+2}$, given by $\{x\in\mathbb{R}^{m+2}:\langle x,\mathbf{u}\rangle=0\}$. In fact
\begin{equation}
\langle\alpha,\mathbf{u}\rangle'=\langle\mathbf{t},\mathbf{u}\rangle=0\Rightarrow \langle\alpha,\mathbf{u}\rangle=c\,\,\mbox{ constant}\,.
\end{equation}
The constant $c$ must be zero. Otherwise, $\alpha$ would be contained on an intersection of $\mathbb{S}^{m+1}(r)$, or $\mathbb{H}^{m+1}(r)$, with a hyperplane not passing through the origin, which is a geodesic sphere \cite{Spivak1979v4}. Since the normal development of a spherical curve does not pass through the origin, we conclude that $c=0$.
\qed

\section{Concluding remarks}

In this work, we furnished necessary and sufficient conditions for a curve to lie on the hypersurface of a geodesic sphere or totally geodesic hypersurface on a hyperbolic space or on a sphere by means of rotation minimization frames. It would be desirable to extend our investigations to the more general setting of Riemannian manifolds that are not necessarily of constant curvature. In this context,  the important concept of normal curves is only valid locally, i.e., one must take into account the injectivity radius of the corresponding exponential map. In addition, it is worth mentioning that a Frenet-like theorem, i.e., two curves are congruent if and only if they have the same curvatures, is valid only for manifolds of constant curvature  \cite{CastrillonLopezAM2015,CastrillonLopezDGA2014}. This may lead to problems in obtaining similar results to ours in terms of the curvatures associated with a given rotation minimizing frame. This is presently under investigation by the authors for some homogeneous spaces and will be the subject of a follow-up work.

\begin{acknowledgements}
The authors would like to thank Gilson S. Ferreira-J\'unior and Gabriel G. Carvalho for useful discussions, the anonymous Referees for their suggestions which have improved the quality of the text, and also the financial support provided by Conselho Nacional de Desenvolvimento Cient\'ifico e Tecnol\'ogico-CNPq (Brazilian agency).
\end{acknowledgements}



\end{document}